\documentclass[10pt]{elsart}
\usepackage{amssymb,amsmath}%,latexsym,amsfonts,amsmath,amsthm}%,verbatim%,hyperref

\newcommand{\newatop}[2]{\genfrac{}{}{0pt}{}{#1}{#2}}
\newcommand{\epf}{$\hfill\qed$}
\newcommand{\pepf}{$\hfill\diamondsuit$}
\newcommand{\eepf}{\hfill\qed}

\newcommand{\lskew}{\!{<}\hspace{-.55em}(\hspace{.2em}}
\newcommand{\rskew}{{>}\hspace{-.65em})\hspace{.2em}}
\newcommand{\func}[2]{#1 \! \left(#2\right)}
\newcommand{\ff}{f\!f}

\def\qplucker{For completeness, it should be noted that the proof of this fact which appears in \cite{TafTow:1} uses a $q$-Laplace expansion in much the same spirit as the classic Laplace expansion was used in the discussion following equation (\ref{eq:grassrels}).}

\def\expansions{This formula is not unique to quantum determinants. Many of the famous noncommutative determinants exhibit this property in some form or another (cf. \cite{GGRW:1}).}

\begin{document}

\begin{frontmatter}
\title{Quantum- and quasi-Pl\"ucker coordinates}
\author{Aaron Lauve\thanksref{nowat}}
\address{Department of Mathematics\\
         Rutgers, The State University of New Jersey\\
         110 Fre\-ling\-huysen Road\\
         Piscatway NJ, 08854-8019, USA}

\thanks[nowat]{Present address: %
LaCIM, 
Universit\'e du Qu\'ebec \`a Montr\'eal,
Case Postale 8888, succursale Centre-ville,
Montr\'eal (Qu\'ebec) H3C 3P8, Canada. \\\texttt{lauve\-@\-lacim.uqam.ca}}

\date{16 November 2005}

\begin{abstract}
We demonstrate a passage from the ``quasi-Pl\"ucker coordinates'' of Gelfand and Retakh, to the quantum Pl\"ucker coordinates built from $q$-generic matrices. In the process, we rediscover the defining relations of the quantum Grassmannian of Taft and Towber and provide that algebra with more concrete geometric origins. 
\end{abstract}

\end{frontmatter}

%-----------------------------------------------------------------------------
% Begin DOCUMENT
%-----------------------------------------------------------------------------
%-----------------------------------------------------------------------------
% Begin SECTION
%-----------------------------------------------------------------------------
\section*{Introduction}
Since the problem of constructing quantum flag and Grassmann spaces was first posed in Manin's Montr\'eal lectures \cite{Man:1}, numerous approaches to the problem have appeared. In this paper, we focus on the efforts of Lakshmibai-Reshetikhin \cite{LakRes:1} and Taft-Towber \cite{TafTow:1} to build the quantized homogeneous coordinate ring $\func{\mathcal{G}_q}{d,n}$ of the Grassmannian of $d$-dimensional subspaces in $K^n$. The difficulty lies in attaching good geometric data to any algebraic structure proposed. 

In this paper we provide further geometric motivation for their definition(s) via the Gelfand-Retakh theory of quasideterminants \cite{GelRet:1}. In 1997, I. Gelfand and V. Retakh introduced coordinates for Grassmannians over division rings in the hope that \emph{specializations} could provide a universal approach to several well-known results in noncommutative geometry. This paper realizes that goal for the quantum Grassmannian of Taft and Towber. We interpret our results as evidence that the definitions of quasi-Pl\"ucker coordinates are the right ones to provide a noncommutative coordinate geometry---and by extention the ``correct'' noncommutative algebra---for many noncommutative settings of interest, not just the quantum setting. 

This paper begins with a review of the classic Grassmannian and its coordinate algebra. We focus our attention on its description in terms of Pl\"ucker coordinates $\{p_I\}$, and Pl\"ucker relations. For example, one has the celebrated identity for minors of a $4\times 2$ matrix $A$:
$$
p_{12}p_{34} - p_{13}p_{24} + p_{23}p_{14} = 0,
$$
where $p_{ij}$ represents the determinant of the submatrix of $A$ formed by taking rows $i$ and $j$ and columns $1$ and $2$. 

The intermediate sections introduce quantum and totally noncommutative versions of this story, the latter relying on \emph{quasi-Pl\"ucker coordinates.} These are certain ratios of quasideterminants which specialize to ratios of minors in the commutative and quantum settings. 

In the final section, we show that the important relations holding among the Pl\"ucker coordinates in the classic and quantum setting are consequences of assorted quasideterminantal identities. For example, if we begin with a ``generic'' $4\times 2$ matrix $A$ and are told that its entries commute with one another, then the identity $(\mathcal P_{1,\{2,3\},\{4\}})$ defined in Section \ref{sec:generic-setting} reduces to
$$
1 = p_{12}p_{32}^{-1}p_{34}p_{14}^{-1} + p_{13}p_{23}^{-1}p_{24}p_{14}^{-1}.
$$
\noindent{\textit{Remark.}} The reader may wish to take a moment to show that the two equations displayed above are equivalent (assuming all symbols $p_{ij}$ are invertible, and $p_{ji}=-p_{ij}$), as it will make some calculations in the sequel more transparent. 

In \cite{LakRes:1} and \cite{LakRes:2} Lakshmibai and Reshetikhin recall the classic realization of $\func{\mathcal G}{d,n}$ as a subalgebra (generated by $d$-minors) of the coordinate algebra for $\mathrm{SL}_n$. With the quantized $\mathcal O_{\mathrm{SL}_n}$ and the quantum determinant provided in \cite{FadResTak:1} on hand, the construction of this algebra is straightforward; cf. \cite{LecZel:1}, \cite{Fio:1},\cite{KelLenRig:1} for modern explorations of its structure. Geometric data appears in the form of representations: they produce its simple modules from the representation theory of $U_q(\mathfrak{sl}_n)$ and use them (along with a modification of Hodge's ``standard monomial theory'' \cite{Hod:1}) to provide a basis for $\func{\mathcal G_q}{d,n}$. 

Taft and Towber \cite{TafTow:1} take a more constructive approach. Beginning with a presentation of $\func{\mathcal G}{d,n}$ by generators and relations, the task was simply to ``quantize'' this presentation to produce the coordinate ring of a quantum Grassmannian. The geometric data here is also indirect: following the suggestion of Faddeev, Reshetikhin, and Takhtajan in \cite{FadResTak:1} they verify their algebra is a comodule algebra over the Hopf algebra $\mathcal O_{\mathrm{SL}_q(n)}$ just as $\func{\mathcal G}{d,n}$ is over $\mathcal O_{\mathrm{SL}_n}$.
They go on to prove that this algebra is the same as the quantum coordinate ring of Lakshmibai and Reshetikhin, strong evidence that indeed this is the ``correct'' quantum $\func{\mathcal G}{d,n}$. 

The aim of this paper is to give more evidence by realizing the generators and relations of Taft and Towber through more geometric considerations. To this end we use quasideterminants. 
Other means of attaching geometric data may be found in \cite{Ohn:1}, where Ohn follows the Artin-Tate-van den Bergh approach to noncommutative projective geometry, and in \cite{Sko:1}, where \v Skoda uses quasideterminant-theory to provide localizations of the quantum algebras in question.

We fix some notation for the remainder of the paper: 
\begin{itemize}
\item[] Fix once and for all, positive integers $d$ and $n$ satisfying $d<n.$ 
\item[] By $[n]$ we mean the set $\{1,2,\ldots n\}$. By $[n]^d$ we mean the set of all $d$-tuples chosen from $[n]$; while ${\scriptstyle \binom{[n]}{d}}$ denotes the set of all subsets of $[n]$ of size $d$. 
\item[] For two integers $n,m$ and two subsets
$I\subseteq[n]$ and $J\subseteq[m]$ we define two common matrices
associated to an $n\times m$ matrix $A$: by $A^{I,J}$ we mean the
matrix obtained by deleting rows $I$ and columns $J$ from $A$; by
$A_{I,J}$ we mean the matrix obtained by keeping only rows $I$ and
columns $J$ of $A$. It will be necessary to simplify the above notation in certain cases: when $I=\{i\}$ and $J=\{j\}$, write $A^{ij}$ in place of $A^{I,J}$; when $|I|=d$ and $J=[d]$, write $A_I$ in place of $A_{I,[d]}$.
\item[] Given two sets $I,J\subseteq[n]$ with $|I|=d, |J|=e$, write $I | J$ for the tuple $(i_1, \ldots, i_d, j_1,\ldots, j_e)$.
\item[] For $\sigma\in\mathfrak{S}_m$, let $\ell(\sigma) = \ell(\sigma1,\sigma2,\cdots,\sigma m)$ denote the \emph{length} of the permutation, i.e. the minimal number of adjacent
swaps necessary to move $(\sigma1,\sigma2,\cdots,\sigma m)$ into $(1,2,\ldots,m)$. Extend $\ell(\cdot)$ to elements of $[n]^m$ in the obvious way; we will make frequent use of $\ell(I\setminus \Lambda | \Lambda)$.
\item[] By $K_q$ we mean an infinite commutative field $K$ of charasteristic $0$ with a distinquished element $q\neq 0$ and $q$ not a root of unity.
\end{itemize}

%-----------------------------------------------------------------------------
% Begin SECTION
%-----------------------------------------------------------------------------
\section{Review of Classical Setting}\label{sec:classic-setting}

\subsection{Determinants}
In this section we work over $\mathbb R$ (cf. \cite{Tow:1} for a treatment over any commutative ring of characteristic $p$ not dividing $d!$). The determinant of a square matrix $A$ will be a main organizing tool in what follows. In addition to the well-known alternating property, the determinant has another property the reader should be familiar with: 

\begin{prop}[Laplace's Expansion] Let $A=\left(a_{ij}\right)_{1\leq i,j\leq m}$. Suppose that $p,p'$ are fixed positive integers with $p+p'=m$, and that  $J=(j_1,\ldots,j_m)$ is a fixed derangement of the columns of $A$. Then 
$$
\big|A\big| = (-1)^{\ell(J)}\sum_{} (-1)^{-\ell(i_1\cdots i_pi'_1\cdots i'_{p'})} \big|A_{\{i_1,\ldots, i_p\},\{j_1,\ldots, j_p\}}\big| \cdot \big|A_{\{i'_1,\ldots,i'_{p'}\},\{j_{p+1},\ldots,j_m\}}\big|
$$
where the sum is over all partitions of $[m]$ into two increasing sets $i_1<\cdots<i_p$ and $i'_1<\cdots<i'_{p'}$. 
\end{prop}
Typically we take $(j_1,\ldots,j_m) = (1,\ldots,m)$, so what's written above is the expansion of the determinant down the first $p$ columns of $A$.

\subsection{Grassmannian}\label{sec:real-flags}
First we recall the embedding of the Grassmannian $Gr(d,n)$ into $\mathbb
P^{{\scriptstyle \binom{n}{d}}-1}$, whose coordinates we will index by
the $d$-subsets of $[n]$. Following \cite{TafTow:1}, we carry
out the construction in $V=(\mathbb R^n)^*$, not in $\mathbb R^n$.

Given a basis $\mathfrak B=\{f_1,\ldots f_n\}$ for $V=\mathbb R^n$, we
will represent a vector $v\in V^*$ as a $n$-tuple
$\left(v_1,\ldots,v_n\right)^{\scriptscriptstyle T}$ where $\langle v,f_i\rangle = v_i$. Any 
$d$-plane $\Gamma\in Gr(d,n)$ can be represented by any $d$
linearly independent vectors within $\Gamma$. We may arrange them as
columns in an $n\times d$ matrix via the coordinatization above. It is clear that any two such matrices $A,B$ represent the same $\Gamma$ if and only if there is an element $g\in\mathrm{GL}_d(\mathbb R)$ satisfying $A = B\cdot g$.

One next forms the map $\eta:Gr(d,n)\rightarrow \mathbb P(\mathbb
R^{{\scriptstyle \binom{n}{d}}})$ as follows. For each $\Gamma$, take
any matrix representation $A$ and map it to the $\binom{n}{d}$-tuple
of its maximal minors. If $A$ and $B$ as above represent the same $\Gamma$, their images will differ only by the scalar $\det g$. Moreover, a matrix $A$ represents an element
of $Gr(d,n)$ if and only if at least one maximal minor is nonzero. One
concludes that $\eta$ is well-defined and injective. (This is the Pl\"ucker 
embedding, and we call the coordinates of $p=\left(p_{\{1,\ldots,d\}}:\cdots:p_{\{n-d+1,\ldots,n\}}\right)\in \mathbb P^{{\scriptstyle \binom{n}{d}}-1}$ the \emph{Pl\"ucker coordinates}.) 

\begin{prop} A point $p\in \mathbb P^{{\scriptstyle \binom{n}{d}}-1}$ belongs to the image of ${\eta}$ if and only if for all $1\leq r\leq d$ and all choices $I\in\binom{[n]}{d+r}$, $J\in\binom{[n]}{d-r}$, its coordinates satisfy
\begin{equation}\label{eq:grassrels}
0 = \sum_{\newatop{\Lambda\subseteq I}{|\Lambda|=r}}(-1)^{\ell(I\setminus \Lambda | \Lambda)}p_{i_1\cdots\hat{i}_{\lambda_1}\cdots\hat{i}_{\lambda_r}\cdots
i_{d+r}}p_{i_{\lambda_1}\cdots i_{\lambda_r} j_{1}\cdots j_{d-r}}\,.
\end{equation} 
\end{prop}
These relations take on many equivalent forms, but as written, they shall be called the \emph{Young symmetry relations} $(\mathcal Y_{I,J})_{(r)}$. The reader may find a proof in \cite{HodPed:1}, one component of which is the ``Basis Theorem'' below. Another component is revealed upon inspection of the following determinant:
$$\left|\begin{array}{cccccc}
a_{i_1,1} & \cdots & a_{i_1,d} & a_{i_1,1} & \cdots & a_{i_1,d} \\
\vdots & & \vdots & \vdots & &\vdots \\
a_{i_{d+r},1} & \cdots & a_{i_{d+r},d} & a_{i_{d+r},d} & \cdots &
a_{i_{d+r},d} \\
0 & \cdots & 0 & a_{j_1,1} & \cdots & a_{j_1,d}\\
\vdots & & \vdots & \vdots & &\vdots \\
0 & \cdots & 0 & a_{j_{d-r},1} & \cdots & a_{j_{d-r},d}
\end{array}\right|.
$$
\noindent{\textit{Remark.}} (a) Use a Laplace expansion down the first $d$ columns to see that this determinant takes the form of (\ref{eq:grassrels}). (b) Subtract the top-left block from the top-right block and discover a hollow matrix, i.e., this determinant is zero.

\subsection{Coordinate Algebra}
There is one technical detail left unsaid after (\ref{eq:grassrels}). In the case $I\cap J\neq\emptyset$, the expressions $p_{i_{\lambda_1}\cdots i_{\lambda_r}j_1\cdots j_{d-r}}$, will not all correspond to subsets of $[n]$. Moreover, order is important. We need $p_{\{i,j,\ldots\}} = -p_{\{j,i,\ldots\}}$, etc. We extend the coordinate functions $\{f_I\}$ on the Pl\"ucker coordinates to $\left\{f_I \mid I\in [n]^d\right\}$ and add the alternating relations $(\mathcal A_I):$
$$
f_I = \left\{\begin{array}{ll}
\mathrm{sgn}(\sigma) f_J & \hbox{ if }\sigma(I)=J\\
0 &\hbox{ if two indices are identical.}
\end{array}\right.
$$
We are now ready to make the

\begin{defn} The homogeneous coordinate ring of $\func{\mathcal G}{d,n}$ is the quotient algebra $\mathbb R\left[f_I\mid I\in {\scriptstyle {[n]}^{d}}\right]/(\mathcal A_{I}; \mathcal Y_{I,J})$.
\end{defn}
\noindent The following theorem suggests we needn't quotient out by a larger ideal.

\begin{thm}[Basis Theorem \cite{HodPed:1}]\label{basis-thm}
If $F$ is any homogeneous polynomial in $f_I$ (modulo $(\mathcal A_{I})$) such that $F(p)=0$ for all $p=p(\Gamma)\in Gr(d,n)$, then $F$ is algebraically dependent on the Young symmetry relations; i.e., 
$$F(p)=\sum_{\newatop{I,|I|=d+1}{J,|J|=d-1}} H_{I,J}(p)\cdot Y_{I,J}(p)\,,$$
where $Y_{I,J}$ is the homogeneous expression appearing on the right-hand side in $(\mathcal Y_{I,J})_{(1)}$ and $H_{I,J}$ is a homogeneous polynomial in the coordinate functions $f_{I}$.
\end{thm}

Note that, interpreting $f_I$ as $\det(A_I)$, we have that any homogeneous polynomial $F$ of degree $m$ in the $f_I$ satisfies $F(A\cdot g) = F(A)(\det q)^m$ as we expect. In the coming sections, we will mimic the constructions above as best as possible.

%-----------------------------------------------------------------------------
% Begin SECTION
%-----------------------------------------------------------------------------
\section{Quantum Setting}\label{sec:quantum-setting}

\subsection{Quantum Determinants}\label{sec:q-det-properties}
Before we introduce the $q$-deformed version of the picture above, we
recall several facts about quantum matrices and quantum determinants. The
reader may find verification of all unproven statements within this
section in \cite{FadResTak:1}, \cite{TafTow:1}, or \cite{Tak:1}.

\begin{defn} An $n\times m$ matrix $X=(x_{ij})$ is called
\emph{$q$-generic} if its entries satisfy all possible relations of the four types below:
\begin{eqnarray}\label{eq:q-genmat}
x_{kj}x_{ki} &=& qx_{ki}x_{kj}\quad(i<j)\\
x_{jk}x_{ik} &=& qx_{ik}x_{jk}\quad(i<j)\\
x_{jk}x_{il} &=& x_{il}x_{jk}\quad(i<j;k<l)\\
x_{jl}x_{ik} &=& x_{ik}x_{jl} + \left(q-q^{-1}\right)x_{il}x_{jk}\quad (i<j;k<l).
\end{eqnarray}
\end{defn}
\noindent{\textit{Remark.}} Any submatrix of a $q$-generic matrix is again $q$-generic.

We let $\mathrm M_{n\times m}(q)$ denote the set of all such $X$. It is a subset of the set of all $n\times m$ matrices with entries in $R$---the often unenunciated ring of study.

Recall that in commutative linear algebra, one can build the inverse
of a matrix $A$ using the determinant:
\begin{equation}\label{build-an-inverse}
\left(A^{-1}\right)_{ij} = (\det A)^{-1}(-1)^{j-i}\det A^{ji}.
\end{equation}
The quantum determinant of a matrix $X=(x_{ij})$ is defined so as to
produce the inverse of a $q$-generic matrix in the same fashion.

\begin{defn} For any square matrix $A=(a_{ij})$ of size $n$, the \emph{quantum determinant} $\det_qA=|A|_q$ is defined by
$$\left|A\right|_q = \sum_{\sigma\in S_n}(-q)^{-\ell(\sigma)}a_{1\sigma1}a_{2\sigma2}\cdots a_{n\sigma n}.$$
\end{defn}

\noindent{\textit{Notation.}} For a subset $I$ of size $m$, we will frequently use
$[I]$ to represent $\det_q\left(A_{I,\{1,\ldots,m\}}\right)$ in order
to simplify notation.

\begin{prop}[Properties of quantum matrices]\label{qdet-properties} Let $X=(x_{ij})$ and $Y=(y_{kl})$
be $q$-generic, with $X$ square and $XY$ defined.
\begin{enumerate}
\item\label{central} The element $\det_qX$ is central in the algebra $K_q\left<x_{ij}\right>/\left(\hbox{\textit{q}-generic relations}\right)$.
\item\label{multiplicative} If $X,Y$ additionally satisfy
$x_{ij}y_{kl}=y_{kl}x_{ij}\, \forall i,j,k,l$ then $XY$ is still $q$-generic; moreover, if $Y$ is square, $det_q(XY)=\det_qX\det_qY$. 
\item\label{inverse} The matrix
$S(X):=\left((-q)^{j-i}\det_qX^{ji}\right)$ satisfies
$S(X)\cdot X=X\cdot S(X)=(\det_qX)I_n, \hbox{the identity matrix}$.
\end{enumerate}
\end{prop}
\noindent{\textit{Warning.}} If $X\in\mathrm M_{n\times n}(q)$ then $X^{-1}\not\in\mathrm M_{n\times n}(q)$; rather it is a member of $\mathrm M_{n\times n}(q^{-1})$.  

\noindent{\textit{Remark.}} Item \ref{central} suggests that $(\det_qA_I)(\det_qA_J)=(\det_qA_J)(\det_qA_I)$ whenever $J\subseteq I$. This will be quite useful in the sequel.

For all $1\leq m\leq n$, define $\mathrm{GL}_q(m)$ to be $\mathrm{GL}_m(R) \cap\mathrm M_{m\times m}(q)$---the $q$-generic matrices which are invertible over $R$. There is not a true group or semigroup structure on this set, e.g. if $X$ is $2\times 2$ $q$-generic, then $X^2$ is not. However, {Proposition \ref{qdet-properties}} suggests that a trace of the desired structure remains: $X\cdot Y\in\mathrm{GL}_q(m)$ when the coordinates of $X$ commute with those of $Y$.

There are two more properties of $\det_q$ which we will need. The first is the ``$q$-alternating'' property of Taft and Towber \cite{TafTow:1}.

\begin{thm} Suppose $X$ is an $n\times n$ $q$-generic matrix, and $A$ is built by choosing rows $i_1,\ldots,i_n$ (not necessarily distinct) from $X$. Then 
\begin{equation}\label{eq:q-alternating}
\mathrm{det}_qA = \left\{\begin{array}{ll}
(-q)^{-\ell(i_1\cdots i_n)}\det_qX & \hbox{if all rows are distinct}\\
0 & \hbox{otherwise}
\end{array}\right..
\end{equation}
\end{thm}
The second property is that often two quantum minors ``$q$-commute:'' 
\begin{defn} Two quantum minors $[I]$ and $[J]$ of a $q$-generic matrix $X$
are said to \emph{$q$-commute} if there is an integer $b$ so that $[J][I] = q^b[I][J]$.
\end{defn}
\noindent For example, we have this
\begin{prop}\label{q-commute-weak} Suppose $i,j\in[n]$ and $M\subset[n]$, with $|M|<n$ and $i<j$. Then the quantum minors $[i\cup M]$ and $[j\cup M]$ satisfy
\begin{equation}\label{eq:q-commute-weak}
[j\cup M][i\cup M] = q [i\cup M][j\cup M].
\end{equation}
\end{prop}
LeClerc and Zelevensky actually prove a much stronger result in \cite{LecZel:1}---giving necessary and sufficient conditions on subsets $I,J$ in order that $[I]$ and $[J]$ $q$-commute. However, their proof involves machinery from \cite{TafTow:1} which we wish to avoid. We present a simple proof of this weak-$q$-commuting property in {Section \ref{sec:quasi-to-quantum}}.

\subsection{Quantum Space}
We are now ready to $q$-deform the picture in {Section \ref{sec:real-flags}.} We move from a vector space over $\mathbb R$ to $n$-dimensional ``quantum space'' $V_q$ over the field $K_q$. We begin by considering a vector space $D^n$ with basis $\mathfrak B=\{f_1,\ldots,f_n\}$, where $D$ is some (unspecified) division algebra over $K_q$. We take $V$ as the left $D$-vector space $V=(D^n)^*=\mathrm{Hom}_D(D^n,D)$; again we build coordinates for vectors $v\in V$ from their behavior on $\mathfrak B$.

We will call a \textbf{point} in $V$ \emph{$q$-generic} if its coordinates satisfy $v_jv_i=qv_iv_j\,(\forall j>i)$. These are the points we wish to study; we call this set $V_q$. {\sc Warning:} this is \emph{not} a vector space over $K_q$ (or $D$) as it is not closed under addition. However, the $K_q$-action inherited from $D$ (it being a $K_q$-algebra) is well-defined. For if  $\alpha\in K_q$, and $v=(v_1,\ldots,v_n)^{\scriptscriptstyle T}\in V_q$, then $\alpha\cdot v\in V_q$ as well (e.g. $(\alpha v_2)(\alpha v_1) =\alpha v_2 v_1\alpha = \alpha qv_1 v_2\alpha = q(\alpha v_1)(\alpha v_2)$).

We will call a $d$-dimensional \textbf{subspace} $W$ of $V$
\emph{$q$-generic} if there is a linearly independent set $\{v^1,\ldots,v^d\}\in V_q\cap W$ so that $A=\left[v^1|\cdots|v^d\right]\in\mathrm M_{n\times d}(q)$. As in the commutative case, $A$ will represent a point in $Gr_q(d,n)$. 

\subsection{Quantum Grassmannian}\label{sec:quantum-flags}
Finally, we define $Gr_q(d,n)$ as a quotient of $\mathrm M_{n\times d}(q)$. We take $A\sim B$ if there is a finite sequence of matrices $\{X_i\}_{1\leq i\leq t}$ chosen from $\mathrm{GL}_q(d)\cup\mathrm{GL}_{q^{-1}}(d)$ satisfying: (i) $(A\cdot X_{1}\cdots X_{i-1})\cdot X_i$ is $q$-generic for all $1\leq i\leq t$; (ii) $B=AX_1\cdots X_t$; (iii) $\det_q B_I = (\det_q A_I)\cdot(\det_q (X_1\cdots X_t ))$ for all $I\in{\scriptstyle \binom{[n]}{d}}$.

\begin{defn} We have defined $Gr_q(d,n)$ above in terms of matrices. We would like to have a  coordinates version as in the commutative case. We identify a point $\Gamma_q$ in the \emph{quantum Grassmannian} $Gr_q(d,n)$ with the set of maximal quantum minors of $A(\Gamma_q)$---its \emph{quantum Pl\"ucker coordinates}.
\end{defn}

\noindent{\textit{Remark.}} Condition (iii) above is fairly restrictive, but it allows us to safely identify two sets of coordinates up to a scalar. We will see shortly that even this is not restrictive enough to completely mimic the classical setting.

From {Section \ref{sec:q-det-properties}} it is clear that a coordinate $[I]$ of $\Gamma_q$ is $q$-alternating in $I$. The coordinates also satisfy a quantized version of the Young symmetry relations\footnote{\qplucker}.

\subsection{Quantized Coordinate Algebra}
Following the classical picture outlined above, we make the
\begin{defn}\label{def:quant-flag} Put $\func{\mathcal G_q}{d,n}=K_q\left<f_I \mid I\in[n]^d\right>/\left(\mathcal A_I;\,\mathcal Y_{I,J}\right)$ where $(\mathcal A_I)$, and $(\mathcal Y_{I,J})$ are now appropriate quantized versions of those from {Section \ref{sec:real-flags}}:
\smallskip

\noindent the alternating relations $(\mathcal A_I)$
$$
f_I = \left\{\begin{array}{ll}
(-q)^{-\ell(I)} f_{\sigma(I)} & \hbox{ if }\sigma\hbox{ orders the entries of }I \\
0 & \hbox{ if two indices are identical.}
\end{array}\right.
$$
\smallskip

\noindent the Young symmetry relations $(\mathcal Y_{I,J}) $
\begin{equation}\label{eq:q-young-symm}
0 = \sum_{\newatop{\Lambda\subseteq I}{|\Lambda|=r}}(-q)^{-\ell(I\setminus \Lambda | \Lambda)} f_{I \setminus \Lambda}f_{\Lambda | J}
\end{equation}
for all $1\leq r\leq d$, $I\in\binom{[n]}{d+r}$, and $J\in\binom{[n]}{d-r}$.
\end{defn}

In \cite{TafTow:1} Taft and Towber give this same definition for the homogeneous coordinate ring of the quantum Grassmannian. They go on to prove a quantized version of the basis theorem: 
\begin{center}\begin{minipage}{.8\textwidth}
\textit{the subalgebra inside $K_q\left<x_{ij} \mid \hbox{\textit{q}-relations}\right>$ generated by $\{[I]\}_{|I|= d}$ is isomorphic to $\func{\mathcal G_q}{d,n}$.} \end{minipage}\end{center}
So not only are the $f_I$ well-defined functions on the points $\Gamma_q\in Gr_q(d,n)$, it would seem $\func{\mathcal G_q}{d,n}$ is the biggest quotient algebra of $K_q\left<f_I\right>$ with this property. 

\noindent{\textit{Remark.}} Note that when we interpret $f_I$ as $\det_q(A_I)$ we have $f_I(A\cdot g)=f_I(A)\det_q(g)$ whenever $A\cdot g\sim A$. Suppose we additionally know that the entries of $g$ commute with those of $A$, then if $F$ is any homogeneous polynomial in $\func{\mathcal G_q}{d,n}$ or degree $m$ we have $F(A\cdot g)=F(A)\cdot(\det_q(g))^m$. This seems to be as close to the classical case as we can come\ldots and not even this is true if we do not add this assumption about $g$. However, we may make a more satisfactory comparison to the classical case when we consider ``homogeneous degree zero'' rational functions in the $\{f_I\}$ (cf. {Proposition \ref{q-rational-invariance}}).

The algebra $\func{\mathcal G_q}{d,n}$ has been well studied since its introduction (cf. \cite{Fio:1,GooLen:1,KelLenRig:1,TafTow:1}). In this paper we concentrate on $Gr_q(d,n)$ itself.

\subsection{Young Symmetry Relations, Simplified}\label{sec:young-symm}
In the classical construction of $\func{\mathcal G}{d,n}$ it is known that all relations of the type in (\ref{eq:grassrels}) with $r>1$ are direct consequences of those with $r=1$ (cf. \cite{HodPed:1} and \cite{Tow:1}). The proofs published there rely heavily on the commutativity of the Pl\"ucker coordinates $\{p_I\}$. What follows is a proof of the same fact for quantum Pl\"ucker coordinates. In addition to giving a new proof for the classical case (set $q=1$), it represents the key lemma for what follows in \textit{Section \ref{sec:quasi-to-quantum}}. 

\noindent{\textit{Notation.}} Given an ordered set $L$ of size $n$ and its $r$-th element $l_r$, let $L_{(r)}$ denote $L\setminus\{l_r\}$. In the event that $l_r\not\in L$ we interpret $L_{(r)}$ as simply a reminder of this fact (i.e. $L_{(r)} = L$). For two subsets $A=\{a_1,\ldots,a_s\}$ and $B=\{b_1,\dots,b_t\}$ of $\{1,\ldots,n\}$, let $[A|B]$ denote $\det_q T_{\{a_1,\ldots,a_s,b_1,\ldots,b_t\},\{1,\ldots,s+t\}}$ for some $q$-generic matrix $T$.

\begin{prop} Let $I,J$ be ordered subsets of $[n]$ with respective sizes $d+r$ and $d-r\,(1\leq r\leq d\leq n)$. Then $(\mathcal Y_{I,J})_{(r)}$ can be written in terms of relations of type $(\mathcal Y_{L,M})_{(r-1)}$. Specifically,
$$\sum_{s=1}^{d+r}(-q)^{2(r-1)-\ell(I_{(s)}|i_{s})}\sum_{\newatop{\Lambda_{(s)}\subset I_{(s)}}{|\Lambda_{(s)}|=r-1}}(-q)^{-\ell(I_{(s)}\setminus\Lambda_{(s)}|\Lambda_{(s)})}\left[I_{(s)}\setminus\Lambda_{(s)}\right]\left[\Lambda_{(s)}|i_{s}|J\right]$$
$$=\left(\sum_{t=0}^{r-1}(-q)^{2t}\right)\sum_{\newatop{\Lambda\subset I}{|\Lambda|=r}}(-q)^{-\ell(I\setminus\Lambda|\Lambda)}\left[I\setminus\Lambda \right]\left[\Lambda|J\right].$$
\end{prop}
\begin{pf}We simply take an arbitrary $\Lambda$ and compare the coefficients on the left- and right-hand sides of the monomial $\left[I\setminus\Lambda\right]\left[\Lambda|J\right]$.
\smallskip
\par\noindent\textit{left-hand side:}
$$\sum_{i_s\in\Lambda} (-q)^{2(r-1)-\ell(I_{(s)}|i_s)}(-q)^{-\ell(I_{(s)}\setminus\Lambda_{(s)}|\Lambda_{(s)})}\left[I\setminus\Lambda\right]\left[\Lambda_{(s)}|i_s|J\right]$$
$$ =\left( \sum_{i_s\in\Lambda} (-q)^{2(r-1)-\ell(I_{(s)}|i_s)-\ell(I_{(s)}\setminus\Lambda_{(s)}|\Lambda_{(s)})-\ell(\Lambda_{(s)}|i_s)}\right)[I\setminus\Lambda][\Lambda|J]$$
\textit{right-hand side:}
$$\left(\sum_{t=0}^{r-1} (-q)^{2t-\ell(I\setminus\Lambda|\Lambda)}\right)[I\setminus\Lambda][\Lambda|J].
$$
Multiplying both sides by $(-q)^{+\ell(I\setminus\Lambda|\Lambda)}$ and using $\ell(I\setminus\Lambda|\Lambda)=\ell(I\setminus\Lambda|\Lambda_{(s)})+\ell(I_{(s)}|i_s)-\ell(\Lambda_{(s)}|i_s)$, we are left with showing
\begin{equation*}
\sum_{s=0}^{r-1} (-q)^{2(r-1)-2\ell(\Lambda_{(s)}|i_s)} = \sum_{t=0}^{r-1} (-q)^{2t}.
\end{equation*}
But $(r-1)-\ell(\Lambda_{(s)}|i_s)$ is exactly $s$.\eepf
\end{pf}

Repeated application of this reduction proves the following important modification to the quantized basis theorem.
\begin{cor}\label{q-youngsymm-short} Equation (\ref{eq:q-young-symm}) in the definition of the $\func{\mathcal G_q}{d,n}$ can be replaced with an abbreviated version---taking only $r=1$.
\end{cor}

\noindent{\textit{Remark.}} (a) Note that this proof fails to work if $q^2$ is an $r$-th root of unity. In the case $q=1$ it additionally fails if the characteristic of the field is $r$. Thus there is no improvement to the situation addressed in \cite{Tow:1} in the commutative case. (b) The lemma was proven for $(|J|,|I|)=(d-r,d+r)$, but a generalization to the setting $(|J|,|I|)=(s-r,t+r)$ with $0\leq r\leq s\leq t\leq d$ is immediate. This extended identity will be utilized in a later paper when we address noncommutative flags.

%-----------------------------------------------------------------------------
% Begin SECTION
%-----------------------------------------------------------------------------
\section{Generic Setting}\label{sec:generic-setting}

\subsection{Quasideterminants}

Gelfand and Retakh suggest that the quasideterminant should be a main organizing tool in noncommutative mathematics; and indeed it has already provided explicit formulas to a variety of noncommutative problems (finding Casimir elements \cite{GKLLRT:1}, \cite{MolRet:1} and factoring noncommutative polynomials \cite{EtiGelRet:1},\cite{GelRet:4} are two notable examples). The results of this paper provide further support for this suggestion. 

The computations in this subsection will be done in the free skew field $K\lskew a_{ij}\rskew$ (cf. \cite{Coh:1}) built on the entries of a matrix $A$ with distinct noncommuting indeterminants. As the definition will make clear, if we instead work with $A$ over an arbitrary noncommutative ring $R$ some quasideterminants may not be defined. A careful study of \cite{Coh:1} reveals that quasideterminants are elements of certain localizations of $R$. The reader will find a more thorough treatment of the quasideterminant and its properties, including some of the proofs omitted below, in \cite{GGRW:1} and \cite{KroLec:1}.

\begin{defn}[Quasideterminant, I] An $n\times n$ matrix $A$ has in general $n^2$ quasideterminants, one for each position in $A$. The
$(ij)$-quasideterminant is defined as follows:
$$|A|_{ij}=a_{ij}-\sum_{r\neq i,s\neq j}
a_{is}\left(|A^{ij}|_{rs}\right)^{-1}a_{rj}.$$
\end{defn}

One may use this definition and (\ref{build-an-inverse}) to easily conclude that in the commutative case, the quasideterminant specializes to the ratio of two determinants: 
$$|A|_{ij} = (-1)^{i+j}(\det A)/(det A^{ij}).$$

\noindent{\textit{Notation.}} It will be convenient to denote the $(ij)$-quasideterminant in another form:
$$\left|\begin{array}{ccc}
 & \vdots & \\
\cdots & a_{ij} & \cdots \\
 & \vdots &
\end{array}\right|_{ij} =
\left|\begin{array}{ccc}
 & \vdots & \\
\cdots & \framebox[1.2\width]{$a_{ij}$} & \cdots \\
 & \vdots &
\end{array}\right|.$$

There is an alternate definition which we will also have occasion to
use. Let $\xi$ be the $i$-th row of $A$ with the $j$-th coordinate
deleted; and let $\zeta$ be the $j$-th column of $A$ with the $i$-th
coordinate deleted. 

\begin{defn}[Quasideterminant, II] For $A,\xi,\zeta$ as above, the $(ij)$-quasideterminant
is defined as follows:
$$|A|_{ij} = a_{ij} - \xi(A^{ij})^{-1}\zeta.$$
\end{defn}

In attempting to make these two definitions agree, one stumbles upon the
first fundamental fact about quasideterminants, 
\begin{equation}\label{eq:inverse}
\left(|A|_{ij}\right)^{-1} = (A^{-1})_{ji},\end{equation}
when the right-hand side is defined and not equal to zero.

The quasideterminant is extremely well-behaved for being a non-commutative determinant (or rather ratio of two). Consider its behavior under elementary transformations of columns.

\begin{prop}\label{elem-trans} Let $A=(a_{ij})$ be a square matrix.
\begin{itemize}
\item (Column Permutations) Suppose $\tau\in S_n$ and $P_{\tau}$ is
the associated (column) permutation matrix. Then $|AP_{\tau}|_{i,\tau j}=|A|_{i,j}$.
\item (Rescaling Columns) Let $B$ be the matrix obtained from
$A$ by multiplying its $r$th column by $\rho$ on the
right. Then 
$$|B|_{ij} =  \left\{\begin{array}{lll}
|A|_{ij}\,\rho & \quad & \hbox{if }j=r \\
|A|_{ij} & \quad & \hbox{if }j\neq r\hbox{ and }\rho\hbox{ is
invertible}.
\end{array}\right.$$
\item (Adding to Columns) Let $B$ be the matrix obtained from $A$ by
adding column $r$ (multiplied on the right by a scalar $\rho$) to column $s$. Then $|B|_{ij}=|A|_{ij}$ if $j\neq r$.
\end{itemize}
\end{prop}

See \cite{GGRW:1} for more details (and for row versions of all the
properties in this subsection). With these properties, we may easily
deduce the following

\begin{prop}\label{combozero} If $A$ is a square matrix and column
$s$ of $A$ is a right-linear combination of the other columns, then $|A|_{rs}=0$ (whenever it is defined). 
\end{prop}
\noindent{\textit{Remark.}} A row version of this is true as well, and will be used below.

\begin{pf}
Through a sequence of steps $A=A(0),\ldots,A(m)=B$, column-reduce $A$ to a matrix $B$: $\mathrm{col}_s(B)=0;\,\mathrm{col}_j(B) = \mathrm{col}_j(A)\,(j\neq s).$ Then {Proposition \ref{elem-trans}} above indicates
$$|A|_{rs} = \left|A(i)\right|_{rs}\quad(\forall 1\leq i\leq m).
$$
Finally, use the second definition of quasideterminant to conclude that $|B|_{rs}$ is indeed zero.\epf
\end{pf}

\begin{prop}[Column Homological Relations] Let $A=(a_{ij})$ be a
square matrix. Then
$$-|A^{kj}|_{il}^{-1}\cdot|A|_{ij} =
|A^{ij}|_{kl}^{-1}\cdot|A|_{kj}\quad(\forall l\neq j).$$
\end{prop}

We will also find a use for the following identity of Krob and LeClerc,  which gives a \emph{one-column} Laplace expansion of the quasideterminant.

\begin{prop} For $A=(a_{ij})$, the $(ij)$-quasideterminant has the
following expansion:
\begin{equation}\label{eq:colexpand}
|A|_{rs} = a_{rs} - \sum_{i\neq r}|A^{is}|_{rl}\cdot|A^{rs}|_{il}^{-1}\cdot
 a_{is}\quad(\forall l\neq s).
\end{equation}
\end{prop}
\begin{pf} From (\ref{eq:inverse}) and the previous proposition we have
\begin{eqnarray*}
1 & = & \sum_{i=1}^{n}|A|_{is}^{-1}\cdot a_{is} \\
|A|_{rs} & = & a_{rs} + \sum_{i\neq r}|A|_{rs}\cdot|A|_{is}^{-1}\cdot
a_{is} \\
|A|_{rs} & = & a_{rs} - \sum_{i\neq
r}|A^{is}|_{rl}\cdot|A^{rs}|_{il}^{-1}\cdot a_{is}. \eepf
\end{eqnarray*}
\end{pf}

\subsection{Noncommutative Pl\"ucker Coordinates}
We may use the quasideterminant to build noncommutative Pl\"ucker
coordinates. One cannot simply replace the determinants appearing
earlier with quasideterminants, because the latter are not invariant
(up to scalar) under $\mathrm{GL}_d$ action. In \cite{GelRet:1,GelRet:3}, Gelfand and Retakh give evidence that certain ratios of quasideterminants are the proper substitute.

\begin{defn}[Quasi-Pl\"ucker Coordinates] Let $A$ be a matrix of size $n\times d\, (n\geq d)$. Let $M$ be a subset of $[n]$ of cardinality $d-1$, and suppose $i,j\in[n]$ with $i\not\in M$. A 
(right-) \emph{quasi-Pl\"ucker coordinate} for $A$ will be defined as $r_{ji}^M(A):=|A_{j\cup
M}|_{js}\cdot|A_{i\cup M}|_{is}^{-1}$ (for any $1\leq s\leq d$).
\end{defn}

\begin{prop}[Compelling Properties] For $A, M, i$,
and $j$, as above, the quasi-Pl\"ucker coordinates satisfy the following:
\begin{itemize}
\item $r_{ji}^M(A)$ does not depend on $s$
\item $r_{ji}^M(A\cdot g)=r_{ji}^M(A)$ for any $g\in\mathrm{GL}_d$
\end{itemize}
\end{prop}

If we associate a point $\Gamma$ in a noncommutative Grassmannian---i.e. a submodule of $V_D=(D^n)^*$ isomorphic to $D^d$ for some division ring $D$---to an $n\times d$ matrix $A$ in a manner similar to what has come before, we might take the quasi-Pl\"ucker coordinates of $\Gamma$ to be the $n^2\binom{n-1}{d-1}$ ``minors'' $r_{ji}^M$. 

Additional nice properties of the $r_{ji}^M$ are worth mentioning.

\begin{prop}\label{rijM-properties} For $A$, $M$, and $i$ as above the following also hold:
\begin{itemize}
\item $r_{ji}^M(A)$ does not depend on the ordering of $M$
\item $r_{ji}^M(A)=\left\{\begin{array}{ll}
0 & \hbox{if }j\in M\\
1 & \hbox{if }j=i\end{array}\right.$
\item $r_{ji}^M\,r_{il}^M = r_{jl}^M \quad (l\not\in M)$
\item $r_{ij}^{M\cup l}\, r_{jl}^{M\cup i}\, r_{li}^{M\cup j} = -1\quad(j,l\not\in M)$
\end{itemize}
\end{prop}

\subsection{Noncommutative Grassmannian}
The fundamental identity holding among the coordinates appears below. It was first observed in \cite{GelRet:3}. We call this identity the ``quasi-Pl\"ucker relations.'' It will allow us to describe Grassmannians and Grassmann algebras in a manner similar to that used in {Section \ref{sec:quantum-flags}}.

\begin{prop}[Quasi-Pl\"ucker Relations] Let $A$ be an $n\times d$
matrix $(n\geq d)$. Then for all subsets $\{i\}, M=\{m_2,\ldots,m_d\}, L=\{l_1,\ldots,l_d\}$ chosen from $\{1,\ldots,n\}$ with $i\not\in M$, we have
\begin{equation}\label{eq:quasi-flag}
(\mathcal P_{i,L,M}):\quad\sum_{j\in L}r_{ij}^{L\setminus j}(A)\cdot r_{ji}^M(A) = 1\,.
\end{equation}
 \end{prop}

\begin{pf}
Using the definition of the quasi-Pl\"ucker coordinates, we show that 
\begin{equation*}
1 = \sum_{j\in L} |A_{i\cup(L\setminus
j)}|_{ir}\cdot|A_{j\cup(L\setminus j)}|_{jr}^{-1}\cdot|A_{j\cup
M}|_{js}\cdot|A_{i\cup M}|_{is}^{-1}\quad (\forall 1\leq r,s\leq d).
\end{equation*}

Let $\xi=\left(\xi_0,\ldots,\xi_t\right)^{\scriptscriptstyle T}$ be
the column vector defined as follows: 
\begin{equation*}
\xi_j = \left\{\begin{array}{ll}
|A_{i\cup M}|_{ir} & \hbox{if }j=0\\
|A_{l_j\cup M}|_{l_jr} & \hbox{otherwise}
\end{array}\right..\end{equation*}
Let $B$ be the matrix $A_{\{i\cup L\},\{1,\ldots,d\}}$ and form the augmented matrix $C = [\xi\vert B]$.

\begin{lem} The matrix $C$ is non-invertible, in particular
$|C|_{11}=0$.
\end{lem}
Using the second definition of quasideterminants, we first notice that
\begin{eqnarray*}
\xi_0 & = & {\left|\begin{array}{ccccc}
a_{i1} & \cdots & \framebox[1.2\width]{$a_{ir}$} & \cdots & a_{id}\\
a_{m_21} & \cdots & a_{m_2r} & \cdots & a_{m_2d}\\
\vdots && \vdots && \vdots\\
a_{m_d1} & \cdots & a_{m_dr} & \cdots & a_{m_dd}
\end{array}\right|}\\
 & = & a_{ir} - \sum_{s\neq r}a_{is}\sum_{t=2}^d \left|
(A_{i\cup M})^{ir}\right|_{m_t s}^{-1}\cdot a_{m_t s}.
\end{eqnarray*}
Computing all of its coordinates at once, we have
\begin{eqnarray*}
\xi &=& \mathrm{col}_r(B) - \mathrm{col}_1(B)\cdot\sum_{t=2}^d \left|
(A_{i\cup M})^{ir}\right|_{m_t 1}^{-1}\cdot a_{m_t 1} - \cdots \\
& & -\,\mathrm{col}_d(B)\cdot\sum_{t=2}^d \left|
(A_{i\cup M})^{ir}\right|_{m_t d}^{-1}\cdot a_{m_t d} \\
 &=& \sum_{j=1}^{d}\mathrm{col}_{j}(B)\cdot\lambda_{j}.
\end{eqnarray*}
Hence the first column is a right-linear combination of the latter columns. In particular, {Proposition \ref{combozero}} implies that $|C|_{11}=0$. \pepf

We next employ (\ref{eq:colexpand}) to $|C|_{11}$ to get the final result: 
\begin{eqnarray*}
0 & = & \xi_0 - \sum_{j=1}^d
|A_{i\cup(L\setminus j)}|_{ir}\cdot|A_{j\cup(L\setminus j)}|_{j
r}^{-1}\cdot\xi_{j}\quad (\forall r)\\
1 & = & \sum_{j= 1}^d
|A_{i\cup(L\setminus j)}|_{ir}\cdot|A_{j\cup(L\setminus j)}|_{jr}^{-1}\cdot|A_{j\cup M}|_{js}\cdot|A_{i\cup M}|_{is}^{-1}\quad (\forall r,s)\\
1 & = & \sum_{j\in L} r_{ij}^{L\setminus j}\cdot r_{ji}^M\,. \eepf
\end{eqnarray*}
\end{pf}

\noindent{\textit{Remark.}} The proof appearing above is new and has an obvious generalization: we only need  $0\leq |M|\leq |L|-1\leq d$ to make the proof work. We will explore this extended identity in a later paper when we address noncommutative flag coordinates. We identify a point $\Gamma$ in the Grassmannian with its collection of quasi-Pl\"ucker coordinates $\big\{r_{ji}^M\big\}$. 

\subsection{Toward a Coordinate Algebra}
One would like a definition of the following sort: the homogeneous coordinate ring of the Grassmannian in the noncommutative setting is the algebra with generators $r_{ij}^M$ and relations all those described above in {Proposition \ref{rijM-properties}} and (\ref{eq:quasi-flag}). However, as all of the symbols are invertible, it seems an algebra of rational functions is more appropriate. In this setting, we have the following theorem.

\begin{thm}\label{g-invt-fcns}
Let $A=(a_{ij})$ be a $n\times d$ matrix with formal entries and let $f(a_{ij})$ be  a rational function over the free skew-field $D$ generated by the $a_{ij}$. Suppose $f$ is invariant under all invertible transformations $A\mapsto A\cdot g,\,(g\in\mathrm{GL}_d(D))$. Then $f$ is a rational function of the quasi-Pl\"ucker coordinates $r_{ij}^M(A)$.
\end{thm}

\begin{pf}
Let $B=A_{\{1,\ldots,d\},\{1,\ldots,d\}}$ and consider the matrix $C=A\cdot B^{-1}$. Then $f(A)=f(C)$, and Gelfand and Retakh have shown that 
\begin{equation*}
(C)_{ij}=\left\{\begin{array}{ll}
\delta_{ij} & j\leq d \\
r_{ij}^{\{1,\ldots,\hat i,\ldots,d\}}(A) & j>d
\end{array}\right.. \eepf
\end{equation*}
\end{pf}

Finally, we would like a version of the basis theorem to be true, e.g. \emph{if $f$ is a rational function in the coordinates $r_{ij}^M$ with $f(A)=0$, then $f=f(\mathcal P_{i,L,M})$ is zero because it can be written in terms of the quasi-Pl\"ucker relations.}  This may be true, but any such theorem is still pending.

%-----------------------------------------------------------------------------
% Begin SECTION
%-----------------------------------------------------------------------------
\section{Quasi $\leadsto$ Quantum}\label{sec:quasi-to-quantum}
In this final section, we return our focus to quantum things (similar results being obtainable for the commutative case via further specialization $q\rightarrow 1$).  

\subsection{Coordinates}

Given a $q$-generic matrix $X$, we have seen that the $(ij)$-th entry of $X^{-1}$ is $(\det_qX)^{-1}(-q)^{j-i}\det_q(X^{ji})$. We have also related the $(ji)$-th quasideterminant of $X$ to the $(ij)$-th entry of
$X^{-1}$. A brief study of this relation yields the following
essential formula\footnote{\expansions}, first introduced in \cite{GelRet:1}:
\begin{equation}\label{eq:quant-via-quasi}
\mathrm{det}_qX = (-q)^{\ell(i_1\cdots i_n)-\ell(j_1\cdots
j_n)}\big|X\big|_{i_1,j_1}\big|X^{i_1j_1}\big|_{i_2,j_2}\cdots|x_{i_nj_n}|_{i_n,j_n};
\end{equation}
moreover, all of the terms on the right-hand side commute with each other {Proposition \ref{qdet-properties}}. We may extend (\ref{eq:quant-via-quasi}) to give quantum determinant expansions for certain matrices associated to $X$.

\begin{prop} Let $A$ be a square matrix, with rows $i_1,\cdots,i_m$ not necessarily ordered
(and not necessarily distinct) chosen from the rows of a $q$-generic matrix $X$. Then 
\begin{equation}\label{eq:quant-via-quasi-ex}
\mathrm{det}_qA = \big|A\big|_{i_11}\big|A^{i_11}\big|_{i_22}\big|A^{i_1i_2,12}\big|_{i_33}\cdots
a_{i_mm} .
\end{equation}
\end{prop}
\begin{pf}If $\det_qA=0$, then the $(i_s)$-th row is the same as some
row $i_t$ ($s<t$) by (\ref{eq:q-alternating})). In this case, $\big|A^{\{i_1\cdots i_{s-1}\},\{i_1\cdots
i_{s-1}\}}\big|_{i_s,i_s}=0$ by {Proposition \ref{combozero}} (row version). 

Otherwise, let $\sigma(j)=i_j$ for $j=1\ldots n$ and use equation
(\ref{eq:q-alternating}) to rewrite (\ref{eq:quant-via-quasi}) as follows:
\begin{eqnarray*}
(-q)^{\ell(\sigma)}\mathrm{det}_qA & = & \mathrm{det}_qX \\
 &=& (-q)^{\ell(\sigma)}\big|X\big|_{\sigma1,1}
\big|X^{\sigma1,1}\big|_{\sigma2,2}\cdots|x_{\sigma n,n}|_{\sigma n,n}\\
 &=& (-q)^{\ell(\sigma)}\big|\sigma^{-1}X\big|_{11}
\big|(\sigma^{-1}X)^{11}\big|_{22}\cdots\big|(\sigma^{-1}X)_{nn}\big|_{nn}\\
 &=& (-q)^{\ell(\sigma)}\big|A\big|_{11}
\big|A^{11}\big|_{22}\cdots|a_{nn}|_{nn}\,,
\end{eqnarray*}
where $\sigma^{-1}$ acts on $X$ by row permutations.
\epf
\end{pf}

\noindent{\textit{Notation.}} For a subset $I$ of size $m$, we will have occasion to use
$\left|i_1\cdots\framebox[1.4\width]{$i_s$}\cdots i_m\right|$ for the
$(i_s1)$-quasideterminant of the matrix $A_{I,\{1,\ldots,m\}}$. For
example, if $B$ is a $2\times 2$ matrix, with rows $i$ and $j$ taken
from some larger matrix $A$, then:
$$|B|_{j1} = \left|\begin{array}{cc}
a_{i1} & a_{i2}\\
\framebox[1.4\width]{$a_{j1}$} & a_{j2}\end{array}\right| = \left|i\,\framebox[1.6\width]{$j$} \,\!\right|.$$

Using this notation---together with the shorthand notation for
$\det_q(A_I)$ described above---the reader may check that the following identities hold:
\begin{itemize}
\item $\left|\framebox[1.4\width]{$i_1$}\cdots i_d\right|=\left[i_1\cdots
i_d\right]\det_q\left(A_{\{i_2,\ldots,i_d\},\{2,\ldots,d\}}\right)^{-1}$
\item $\left|\framebox[1.8\width]{$i$} m_2\cdots m_d\right|\left|\framebox[1.6\width]{$j$} m_2\cdots
m_d\right|^{-1}=\left[i m_2\cdots
m_d\right]\left[j m_2\cdots m_d\right]^{-1}$
\end{itemize}

\begin{pf*}{PROOF, Proposition \ref{q-commute-weak}.} Consider the following column homological relation for the $q$-generic matrix $A_{i\cup j\cup M}$:
\begin{eqnarray*}-\big|A^{j1}\big|_{i2}^{-1}\cdot\big|A\big|_{i1}&=&\big|A^{i1}\big|_{j2}^{-1}\cdot\big|A\big|_{j1}\\
-\big|A^{i1}\big|_{j2}\cdot\big|A^{j1}\big|_{i2}^{-1} &=& \big|A\big|_{j1}\cdot\big|A\big|_{i1}^{-1}.
\end{eqnarray*}
We apply the simple identities above, the $q$-alternating property, and {Proposition \ref{qdet-properties}(\ref{central})} to finish the proof. 
\par\noindent\emph{Left-hand side:}
\begin{eqnarray*}
-\big|A^{i1}\big|_{j2}\cdot\big|A^{j1}\big|_{i2}^{-1} &=&
-\left|\framebox[1.6\width]{$j$}m_2\cdots m_d\right|\left|\framebox[1.6\width]{$i$}m_2\cdots m_d\right|^{-1}\\
&=&-[jm_2\cdots m_d][im_2\cdots m_d]^{-1},
\end{eqnarray*}
using the identities above starting from column 2 of the original matrix $A$.
\par\noindent\emph{Right-hand side:}
\begin{eqnarray*}
\big|A\big|_{j1}\cdot\big|A\big|_{i1}^{-1} &=& \left|i\framebox[1.6\width]{$j$}m_2\cdots m_d\right|\left|\framebox[1.8\width]{$i$}jm_2\cdots m_d\right|^{-1}\\
 &=& \left([jim_2\cdots m_d][im_2\cdots m_d]^{-1}\right)\cdot\left( [ijm_2\cdots m_d][jm_2\cdots m_d]^{-1}\right)^{-1}\\
 &=& [im_2\cdots m_d]^{-1}[jim_2\cdots m_d] \,[ijm_2\cdots m_d]^{-1}[jm_2\cdots m_d]\\
 &=& -q^{\pm1}\,[im_2\cdots m_d]^{-1}[jm_2\cdots m_d],
\end{eqnarray*} 
where the power of $-q$ depends on whether $i<j$ or $i>j$. Note the heavy reliance on the centrality of quantum determinants, {Proposition \ref{qdet-properties}(\ref{central})}. The result now follows by clearing denominators.
\epf
\end{pf*}

\subsection{Grassmannians}
In this section we prove the main result of this paper, the quantized coordinate algebra $\func{\mathcal G_q}{d,n}$ of Taft and Towber results from specializing the geometry of the generic Grassmannian. 

\begin{thm}[Quasi- Specialization]\label{theo:quasi-to-quantum} Let $A$ be an $n\times d$ $q$-generic matrix representing a point $\Gamma_q$ in the quantum Grassmannian $Gr_q(d,n)$. Then all the relations among the coordinates $\left\{[I] \mid I\in{\scriptstyle\binom{[n]}{d}}\right\}$ of $\Gamma_q$ are consequences of the coordinate-relations for its quasi-Pl\"ucker coordinates $\left\{r_{ji}^M(A)\right\}$.
\end{thm}
Alternatively, beginning with an $n\times d$ matrix $A$ of
indeterminants over the free skew field built on the $\{a_{ij}\}$, all
relations of the form ($\mathcal Y_{I,J}^*$) are direct consequences of adding $q$-genericity to the
quasi-Pl\"ucker relations ($\mathcal P_{i,L,M}$) already holding for $A$.

\begin{pf}
Along with $q$-genericity we add its easy consequences---the $q$-alternating and (weak) $q$-commuting properties of equations (\ref{eq:q-alternating}) and (\ref{eq:q-commute-weak}).

We have as our target $\left(\mathcal Y_{IJ}\right)_{(r)}$. By {Corollary \ref{q-youngsymm-short}}, we may assume $r=1$; so let $I=\{i_1,\ldots,i_{d+1}\}$ and $J=\{j_1,\ldots,j_{d-1}\}$.  Starting from the relation $(\mathcal P_{i_1,I_{(1)},J})$ we have:
\begin{eqnarray*}
1 &=& \sum_{j\in L}r_{ij}^{L\setminus j}\cdot r_{ji}^M \\
1 &=& \sum_{j\in L}\big|\framebox[1.8\width]{$i$}L\setminus j\big|\big|\framebox[1.6\width]{$j$}L\setminus j\big|^{-1} \cdot\big|\framebox[1.6\width]{$j$}M\big|\big|\framebox[1.8\width]{$i$}M\big|^{-1}\\
1 &=& \sum_{2\leq\lambda\leq d+1} \big|\framebox[1.4\width]{$i_{1}$}I_{(1)}\setminus i_{\lambda}\big|\big|\framebox[1.2\width]{$i_{\lambda}$}I_{(1)}\setminus i_{\lambda}\big|^{-1} \cdot\big|\framebox[1.2\width]{$i_{\lambda}$}J\big|\big|\framebox[1.4\width]{$i_{1}$}J\big|^{-1}\\
1 &=& \sum_{2\leq\lambda\leq d+1}\big[i_{1}I_{(1)}\setminus i_{\lambda}\big]\big[i_{\lambda}I_{(1)}\setminus i_{\lambda}\big]^{-1} \cdot\big[i_{\lambda}J\big]\big[i_{1}J\big]^{-1}
\end{eqnarray*}
\begin{eqnarray*}
\big[i_{1}J\big] & =& \sum_{2\leq\lambda\leq d+1} {q}\big[i_{\lambda}I_{(1)}\setminus i_{\lambda}\big]^{-1}\big[i_{1}I_{(1)}\setminus i_{\lambda}\big]\big[i_{\lambda}J\big]\\
\big[i_{1}J\big] & =& \sum_{2\leq\lambda\leq d+1} {q}(-q)^{+\ell(i_{\lambda}|I_{(1)}\setminus i_{\lambda})} \big[I_{(1)}\big]^{-1}\big[i_{1}I_{(1)}\setminus i_{\lambda}\big]\big[i_{\lambda}J\big]\\ 
\big[I_{(1)}\big]\big[i_{1}J\big] & =& \sum_{2\leq\lambda\leq t+r} {q}(-q)^{+\ell(i_{\lambda}|I_{(1)}\setminus i_{\lambda})} \big[i_{1}I_{(1)}\setminus i_{\lambda}\big]\big[i_{\lambda}J\big]\\
\big[I_{(1)}\big]\big[i_{1}J\big] & =& -\sum_{2\leq\lambda\leq d+1} (-q)^{+\ell(i_{\lambda}|I\setminus i_{\lambda})}\big[I\setminus i_{\lambda}\big] \big[i_{\lambda}J\big]\end{eqnarray*}
\begin{eqnarray*}
0 &=& (-q)^{-\ell(I_{(1)}|i_{1})}\big[I_{(1)}\big]\big[i_{1}J\big] + \sum_{2\leq\lambda\leq d+1} (-q)^{+\ell(i_{\lambda}|I\setminus i_{\lambda})} (-q)^{-\ell(I_{(1)}|i_{1})}\big[I\setminus i_{\lambda}\big]\big[i_{\lambda}J\big]\\
0 &=&  (-q)^{-\ell(I_{(1)}|i_{1})} \big[I_{(1)}\big]\big[i_{1}J\big] + \sum_{2\leq\lambda\leq d+1}(-q)^{-\ell(I_{(\lambda)}|i_{\lambda})}\big[I_{(\lambda)}\big]\big[i_{\lambda}J\big] \\
0 &=& \sum_{1\leq\lambda\leq d+1}(-q)^{-\ell(I_{(\lambda)}|i_{\lambda})} \big[I_{(\lambda)}\big]\big[i_{\lambda}J\big].
\end{eqnarray*}
This is exactly the targeted $\left(\mathcal Y_{IJ}\right)$. Now, we implicitly began with the assumption $i=i_{1}\not\in J$, but any choice from $I\setminus J$ could have been made for $i$. Finally, if $I\setminus J=\emptyset$, then (\ref{eq:q-young-symm}) reads $0=0$ by the $q$-alternating property.
\epf
\end{pf}

\subsection{Coordinate Algebras}
We conclude this section with the introduction of a natural algebra of functions on the quantum Grassmannian. This algebra is invariant under the relation $\sim$ introduced in {Section \ref{sec:quantum-flags}.} Moreover, its elements $F$ are identically zero on $Gr_q(d,n)$ only if they are zero for quasi-Pl\"ucker reasons. 

In \cite{KelLenRig:1} we learn that $\func{\mathcal G_q}{d,n}$ is a noetherian domain, and as such has a (right, Ore) skew-field of fractions $\mathcal D$. Namely, every element of $\mathcal D$ can be written as $GH^{-1}$ with $G,H\in \func{\mathcal G_q}{d,n}$. This field is too big to be an appropriate field of fractions for $Gr_q(d,n)$; we look for the $\sim$-invariant functions within $\mathcal D$. 

\begin{prop} Let $R$ be a noetherian domain with right field of fractions $D$. If $R$ is graded, then the subset $D_0=\{gh^{-1}\in R \mid g,h \hbox{ are homogeneous of the same degree}\}$ is a well-defined subfield of $D$.
\end{prop}
\begin{pf}
Given $ef^{-1}$ and $gh^{-1}$ in $D_0$, we may add and multiply these two fractions together by the Ore conditions in $D$:
\begin{enumerate}
\item[$(+):$]  We know $\exists u,v\in R$ with $fu = hv$. So we may write $ef^{-1} + gh^{-1}= (eu)(fu)^{-1} + (gv)(hv)^{-1} = (eu + gv)(hv)^{-1}$. 
\item[$(\times):$] we know $\exists u',v'\in R$ with $fu' = gv'$. So we may write $ef^{-1}gh^{-1} = (eu')(fu')^{-1}(gv')(hv')^{-1} = (eu')(hv')^{-1}$.
\end{enumerate}
One question is whether $u,v,u',v'$ may be chosen to be homogeneous elements of $R$. This is straightforward to check:

\par\noindent Write $u = \sum_{i=s}^{\infty} u_i$ and $v=\sum_{j=t}^{\infty} v_j$ (finite sums) with $u_s,v_t\neq 0$ (the pieces of $u$ and $v$ of lowest degree). Now, $fu_s$ is the lowest degree piece of $fu$ because $f$ is homogeneous and $R$ is a (graded) domain. Similarly, $hv_t$ is the lowest degree piece of $hv$. Finally, $fu=hv\Rightarrow fu_s=hv_t$ again by the grading of $R$. Hence we may assume $u$ and $v$ (and $u'$ and $v'$) are homogeneous elements of $R$.

Next, we must ask whether the resulting fractions in $(+)$ and $(\times)$ above belong to $D_0$. Again, this is easy to check, and we do so only for $(+)$:

\par\noindent In the case of $(+)$ we have $\deg e + \deg u = \deg f + \deg u = \deg h + \deg v = \deg g + \deg v$, so $\deg(eu+gv) = \deg(hv)$ as needed.
\epf
\end{pf}

For what remains, we will need a stronger version of the $q$-commuting property than was proved above. Specifically, we need the following identity.

\begin{prop} Put $f_{[-d]}:=f_{\{n-d+1,\ldots, n\}}$. Then for all $I\in{\scriptstyle\binom{[n]}{d}}$, we have 
$$f_{[-d]}f_I  = q^{|[-d]\setminus I|} f_If_{[-d]}.$$\end{prop}
\par\noindent One can find a proof of this well-known identity in \cite{LecZel:1}, which, after the specialization results of the previous section, we are now free to use. 

Define $\mathcal D_0\subset\mathcal D$ as in the proposition. This is the algebra of functions we seek. We must show that: (i) $\mathcal D_0$ is $\sim$-invariant; (ii) it is neither too big nor too small inside $\mathcal D$. 

\par\noindent\textit{Invariance:}\\
Write $\hat f_I$ for $f_If_{[-d]}^{-1}$ inside $\mathcal D$. Note that $\hat f_I$ is $\sim$-invariant. Finally, take $F=GH^{-1}\in\mathcal D_0$ (with $\deg G=\deg H = b$), and write $GH^{-1} = (Gf_{[-d]}^{-b})(Hf_{[-d]}^{-b})^{-1}=\widehat G\widehat H^{-1}$ in $\mathcal D$. Here we have written $\widehat G$ for the rearrangement of $Gf_{[-d]}^{-b}$ putting one factor of $f_{[-d]}^{-1}$ to the right of each symbol $f_I$ appearing in $G$. Then $GH^{-1}(Ag) = \widehat G\widehat H^{-1}(Ag) = \widehat G\widehat H^{-1}(A) = GH^{-1}(A)$ as needed.

\par\noindent\textit{Correct Size:}\\
We look at the fields of fractions on the affine pieces of our projective space $Gr_q(d,n)$. $\mathcal D_0$ should contain them all, and be no bigger than necessary.
Consider the ``affine patch'' of points $X_{[-d]}=\left\{\{|A_{I}|_q\}\,:\, |A_{\{n-d+1,\ldots, n\}}|_q\neq 0 \right\}$ inside $Gr_q(d,n)$; $f_{[-d]}^{-1}$ is a well-defined function here. Moreover, by property (iii) of $\sim$ we have $f_I f_{[-d]}^{-1} (Ag) = f_I f_{[-d]}^{-1} (A)$ when $Ag\sim A$. So we may consider the subalgebra $\mathcal A$ of $\mathcal D$ generated by $f_I f_{[-d]}^{-1}$ as a piece of the field of $\sim$-invariant functions we're looking for. By the previous proposition, we may write every element of $\mathcal A$ as $Gf_{[-d]}^{-\mathrm{deg} G}$ in $\mathcal D$, where $G$ is a homogeneous polynomial in $\func{\mathcal G_q}{d,n}$. Finally, $\mathcal A$ is noetherian (cf. \cite{KelLenRig:1}, {Theorem 1.4}), so we may consider its right field of fractions $\ff\mathcal A\subseteq\mathcal D$. Observe that $\mathcal D_0\subseteq \ff\mathcal A$: given $GH^{-1}\in\mathcal D_0$, we have
$$GH^{-1} = \left({G}{f_{[-d]}^b}\right)\cdot\left({H}{f_{[-d]}^b}\right)^{-1}\in \ff\mathcal A\,.
$$ 
On the other hand, note that all rings corresponding to all affine patches are subalgebras of $\mathcal D_0$, and thus so are their fields of fractions---to whatever extent they exist. So we arrive at the natural 
\begin{defn} The field of functions on $Gr_q(d,n)$ is the subfield $\mathcal D_0$ of $\mathcal D$ generated by all elements $G\cdot H^{-1}$ with $G,H\in\func{\mathcal G_q}{d,n}$ homogeneous of the same degree.
\end{defn}

\begin{prop}\label{q-rational-invariance} If $F\in\mathcal D_0$, then $F$ is a rational function in $\{(f_I)(f_J)^{-1}\,:\, |I\cap J|= d-1\}$. 
\end{prop} 

\noindent{\textit{Remark.}} ``$\sim$'' is too strict a relation to allow Gaussian elimination\ldots a procedure necessary in the proof of  {Proposition \ref{g-invt-fcns}}, so we cannot simply pass from quasi- to quantum- in that proposition.

\begin{pf*}{PROOF, Sketch.}
The proof comes from the special form $F$ takes. Let's consider the commutative case for a moment. Start from $F=G/H$ with $G$ and $H$ homogeneous of the same degree, $b$ say. Here one may divide the top and bottom by $f_{[-d]}^b$ and ``interpolate'' between the coordinate functions $f_I$ occuring in $G$ and $H$ to get this same result in a more elementary fashion.
\begin{exmp}
\begin{eqnarray*}
\frac{f_{\{346\}}+f_{\{123\}}}{f_{\{135\}}} &=& \frac{f_{\{346\}}f_{\{456\}}^{-1}+f_{\{123\}}f_{\{456\}}^{-1}}{f_{\{135\}}f_{\{456\}}^{-1}}\\
&=& \frac{(f_{\{346\}}f_{\{456\}}^{-1})+(f_{\{123\}}f_{\{126\}}^{-1})(f_{\{126\}}f_{\{156\}}^{-1})(f_{\{156\}}f_{\{456\}}^{-1})}{(f_{\{135\}}f_{\{345\}}^{-1})(f_{\{345\}}f_{\{456\}}^{-1})}\,.
\end{eqnarray*} 
\end{exmp}  
In the quantum setting, the same argument works as $f_{[-d]}$ $q$-commutes with every other coordinate function.
\epf
\end{pf*}

We have given some motivation for the further study of $\mathcal D_0$. We conclude this section by showing that, like $Gr_q(d,n)$, it's behavior is governed by its quasi-counterpart. 

\begin{thm} If  $F\in\mathcal D_0$ is identically zero on $Gr_q(d,n)$, then $F$ is zero as a consequence of quasi-Pl\"ucker coordinate considerations.
\end{thm}

\begin{pf}
Let $Y_{I,J}$ denote the right-hand side of (\ref{eq:q-young-symm}) and $P_{i_1,I_{(1)},J}$ denote the left-hand side of (\ref{eq:quasi-flag})---so $Y_{I,J}=0$ in $\func{\mathcal G_q}{d,n}$, and $1-P_{i_1,I_{(1)},J}=0$ in $\mathcal D_0$. For $F\in\mathcal D_0$, write $F=GH^{-1}$ as above, with $G(\Gamma_q)=0,H(\Gamma_q)\neq 0$. Then $G$---by the quantized basis theorem---is in the ideal generated by relations of type $(\mathcal Y_{I,J}^*)_{(1)}$. Write $G$ as such, then consider $\bar G\in\mathcal D_0$ built from $G$ by factoring each expression $w(Y_{I,J}^*)w'$ occuring as $wf_{\{i_1\cdots i_d\}}(1- P_{i,I_{(1)},J})f_{\{i_{d+1}\cdots i_{d+r}j_1\cdots j_{d-r}\}}w'$ in the manner carried out in the proof of {Theorem \ref{theo:quasi-to-quantum}}.
\epf
\end{pf}

%-----------------------------------------------------------------------------
% Begin SECTION
%-----------------------------------------------------------------------------
\section{Future Steps}
As mentioned earlier, we anticipate following this paper with another addressing more general quantum flags. Already from the results of this paper, one may confidently go on to create Grassmannians in other noncommutative settings where amenable determinants exist (e.g. superalgebras). 

Beyond ``specializations'' such as those above, it would be interesting to study the ring of quasi-Pl\"ucker coordinates itself. Recall the classical result: \emph{the homogeneous coordinate ring for the flag variety is a model for the irreducible polynomial representations of $\mathrm{GL}_n$.} One challenge would be to use the quasi-Pl\"ucker coordinates to construct a noncommutative representation theory. %After the evident success of the quasi-Pl\"ucker coordinates, it is clear that this stands as a significant direction for further research.  $K\!\langle r_{ij}^M \mid \hbox{quasi-Pl\"ucker relations }\rangle$ and its localizations have not been addressed in this paper. 

%-----------------------------------------------------------------------------
% Begin SECTION
%-----------------------------------------------------------------------------
\section*{Acknowledgement}
The author would like to thank Vladimir Retakh and Robert Wilson for many helpful discussions, and for encouraging the writing of this paper.

%%%%I used these commands to build the bibliography below
%\bibliographystyle{amsplain}     
%\bibliography{../../global-bib.bib}  

\begin{thebibliography}{10}

\bibitem{Coh:1}
Paul~Moritz Cohn, \emph{Skew field constructions}, London Mathematical Society
  Lecture Note Series, no.~27, Cambridge University Press, Cambridge, 1977.

\bibitem{EtiGelRet:1}
Pavel Etingof, Israel Gelfand, and Vladimir Retakh, \emph{Factorization of
  differential operators, quasideterminants, and nonabelian {T}oda field
  equations}, Math. Res. Lett. \textbf{4} (1997), no.~2-3, 413--425.

\bibitem{Fio:1}
R.~Fioresi, \emph{Quantum deformation of the {G}rassmannian manifold}, J.
  Algebra \textbf{214} (1999), no.~2, 418--447.

\bibitem{GelRet:1}
I.~M. Gel{\cprime}fand and V.~S. Retakh, \emph{Determinants of matrices over
  noncommutative rings}, Funktsional. Anal. i Prilozhen. \textbf{25} (1991),
  no.~2, 13--25, 96.

\bibitem{GelRet:3}
I.~M. Gelfand and V.~S. Retakh, \emph{Quasideterminants, {I}}, Selecta Math.
  (N.S.) \textbf{3} (1997), no.~4, 517--546.

\bibitem{GGRW:1}
Israel Gelfand, Sergei Gelfand, Vladimir Retakh, and Robert~Lee Wilson,
  \emph{Quasideterminants}, Adv. in Math. \textbf{193} (2005), no.~1, 56--141.

\bibitem{GelRet:4}
Israel Gelfand and Vladimir Retakh, \emph{Noncommutative {V}ieta theorem and
  symmetric functions}, The Gelfand Mathematical Seminars, 1993--1995,
  Birkh{\"a}user Boston, Boston, MA, 1996, pp.~93--100.

\bibitem{GKLLRT:1}
Israel~M. Gelfand, Daniel Krob, Alain Lascoux, Bernard Leclerc, Vladimir~S.
  Retakh, and Jean-Yves Thibon, \emph{Noncommutative symmetric functions}, Adv.
  Math. \textbf{112} (1995), no.~2, 218--348.

\bibitem{GooLen:1}
K.~R. Goodearl and T.~H. Lenagan, \emph{Quantum determinantal ideals}, Duke
  Math. J. \textbf{103} (2000), no.~1, 165--190.

\bibitem{Hod:1}
W.~V.~D. Hodge, \emph{Some enumerative results in the theory of forms}, Proc.
  Cambridge Philos. Soc. \textbf{39} (1943), 22--30.

\bibitem{HodPed:1}
W.~V.~D. Hodge and D.~Pedoe, \emph{Methods of algebraic geometry. {V}ol. {I}},
  Cambridge Mathematical Library, Cambridge University Press, Cambridge, 1994,
  Book I: Algebraic preliminaries, Book II: Projective space, Reprint of the
  1947 original.

\bibitem{KelLenRig:1}
A.~C. Kelly, T.~H. Lenagan, and L.~Rigal, \emph{Ring theoretic properties of
  quantum {G}rassmannians}, J. Algebra Appl. \textbf{3} (2004), no.~1, 9--30.

\bibitem{KroLec:1}
Daniel Krob and Bernard Leclerc, \emph{Minor identities for quasi-determinants
  and quantum determinants}, Comm. Math. Phys. \textbf{169} (1995), no.~1,
  1--23.

\bibitem{LakRes:1}
V.~Lakshmibai and N.~Reshetikhin, \emph{Quantum deformations of {${\rm SL}\sb
  n/B$} and its {S}chubert varieties}, Special functions (Okayama, 1990),
  ICM-90 Satell. Conf. Proc., Springer, Tokyo, 1991, pp.~149--168.

\bibitem{LakRes:2}
\bysame, \emph{Quantum flag and {S}chubert schemes}, Deformation theory and
  quantum groups with applications to mathematical physics (Amherst, MA, 1990),
  Contemp. Math., vol. 134, Amer. Math. Soc., Providence, RI, 1992,
  pp.~145--181.

\bibitem{LecZel:1}
Bernard Leclerc and Andrei Zelevinsky, \emph{Quasicommuting families of quantum
  {P}l\"ucker coordinates}, Kirillov's seminar on representation theory, Amer.
  Math. Soc. Transl. Ser. 2, vol. 181, Amer. Math. Soc., Providence, RI, 1998,
  pp.~85--108.

\bibitem{Man:1}
Yu.~I. Manin, \emph{Quantum groups and noncommutative geometry}, Universit\'e
  de Montr\'eal Centre de Recherches Math\'ematiques, Montreal, QC, 1988.

\bibitem{MolRet:1}
Alexander Molev and Vladimir Retakh, \emph{Quasideterminants and {C}asimir
  elements for the general linear {L}ie superalgebra}, Int. Math. Res. Not.
  (2004), no.~13, 611--619.

\bibitem{Ohn:1}
Christian Ohn, \emph{``{C}lassical'' flag varieties for quantum groups: the
  standard quantum {${\rm SL}(n,\bold C)$}}, Adv. Math. \textbf{171} (2002),
  no.~1, 103--138.

\bibitem{FadResTak:1}
N.~Yu. Reshetikhin, L.~A. Takhtadzhyan, and L.~D. Faddeev, \emph{Quantization
  of {L}ie groups and {L}ie algebras}, Algebra i Analiz \textbf{1} (1989),
  no.~1, 178--206, English transl.: \textit{Leningrad Math. J.} \textbf{1}
  (1990), no. 1, 193--225.

\bibitem{Sko:1}
Zoran {\v{S}}koda, \emph{Localizations for construction of quantum coset
  spaces}, Noncommutative geometry and quantum groups (Warsaw, 2001), Banach
  Center Publ., vol.~61, Polish Acad. Sci., Warsaw, 2003, pp.~265--298.

\bibitem{TafTow:1}
Earl Taft and Jacob Towber, \emph{Quantum deformation of flag schemes and
  {G}rassmann schemes, {I}. {A} $q$-deformation of the shape-algebra for
  $\mathrm{{G}{L}}(n)$}, J. Algebra \textbf{142} (1991), no.~1, 1--36.

\bibitem{Tak:1}
Mitsuhiro Takeuchi, \emph{A short course on quantum matrices}, New directions
  in Hopf algebras, Math. Sci. Res. Inst. Publ., vol.~43, Cambridge Univ.
  Press, Cambridge, 2002, Notes taken by Bernd Str\"uber, pp.~383--435.

\bibitem{Tow:1}
Jacob Towber, \emph{Young symmetry, the flag manifold, and representations of
  {${\rm GL}(n)$}}, J. Algebra \textbf{61} (1979), no.~2, 414--462.

\end{thebibliography}
%%%%
\def\cprime{$'$} \def\cprime{$'$} \def\cprime{$'$} \def\cprime{$'$}
\providecommand{\bysame}{\leavevmode\hbox to3em{\hrulefill}\thinspace}
\providecommand{\MR}{\relax\ifhmode\unskip\space\fi MR }
\providecommand{\MRhref}[2]{%
  \href{http://www.ams.org/mathscinet-getitem?mr=#1}{#2}
}
\providecommand{\href}[2]{#2}

\end{document}